\documentclass[12pt]{amsart}
 \usepackage[latin1]{inputenc}
 \usepackage[dvips]{graphicx}
 \usepackage{wrapfig}
 \usepackage{amsmath}
 \usepackage{amsthm}
 \usepackage{amsfonts}
 \usepackage{amssymb}
 \usepackage{layout}
 \usepackage{verbatim}
 \usepackage{alltt}
\usepackage{yfonts}

\usepackage[all]{xy}
\usepackage{setspace}
\newtheorem*{ithm}{Theorem}

\newcommand{\QQ}{\mathbb{Q}}
\newcommand{\FF}{\mathcal{F}}
\newcommand{\lra}{\longrightarrow}
\newcommand{\ZZ}{\mathbb{Z}}
\newcommand{\PP}{\mathcal{P}}
\newcommand{\NN}{\mathcal{N}}
\newcommand{\ra}{\rightarrow}

\newcommand{\be}{\begin{equation}}
\newcommand{\ee}{\end{equation}}

\newcommand{\XX}{\mathcal{X}}
\newcommand{\kk}{\mathcal{K}}
\newcommand{\al}{\mathcal{L}}
\newcommand{\JJ}{\mathbb{J}}

\newcommand{\eee}{\epsilon}
\newcommand{\hli}{\hbox{H}_{\al}}
\newcommand{\hhh}{\bigwedge^{r-1}\textbf{\textup{Hom}}(\mathbb{V},\ZZ_p[[\textup{Gal}(M/k)]])}
\newcommand{\hhhk}{\bigwedge^{r-1}\textbf{\textup{Hom}}(\mathbb{V}^{\chi},\ZZ_p[[\textup{Gal}(M/k)]]^{\chi})}
\newcommand{\qq}{\hbox{\frakfamily q}}
\newcommand{\vv}{\mathbb{V}}
\newcommand{\all}{\mathbb{L}}

\numberwithin{equation}{section}
\newtheorem{thm}{Theorem}[section]
\newtheorem{lemma}[thm]{Lemma}
\newenvironment{define}{\par\medskip\noindent\refstepcounter{thm}
\bgroup{\hspace*{-0.15 cm}\bf{Definition}
\thethm.}\bgroup}{\egroup \egroup\par\medskip}\newtheorem{prop}[thm]{Proposition}
\newtheorem{cor}[thm]{Corollary}
\newenvironment{rem}{\par\medskip\noindent\refstepcounter{thm}
\bgroup{\hspace*{-0.15 cm}\bf{Remark} \thethm.}\bgroup}{\egroup
\egroup\par\medskip} \parskip 2pt

\begin{document}
\title{{K}\lowercase{olyvagin systems of }{S}\lowercase{tark units}}

\author{K\^az\i m B\"uy\"ukboduk}

\address{Kazim Buyukboduk \hfill\break\indent
Department of Mathematics
\hfill\break\indent 450 Serra Mall, Bldg 380 \hfill\break\indent Stanford, CA, 94305
\hfill\break\indent USA} \email{{\tt kazim@math.stanford.edu }\hfill\break\indent {\it
Web page:} {\tt math.stanford.edu/$\sim$kazim}}
\address{Current Address: \hfill\break\indent
IH\'{E}S, Le Bois-Marie. 35,\hfill\break\indent 
Route de Chartres, 
\hfill\break\indent F-91440 Bures-sur-Yvette,
\hfill\break\indent France}
\keywords{Stark conjectures, Euler systems, Kolyvagin systems.}
\subjclass[2000]{11R27, 11R29, 11R42}

\begin{abstract}

In this paper we construct, using Stark elements of Rubin~\cite{ru96}, Kolyvagin systems for certain \emph{modified} Selmer structures (that are adjusted to have core rank \emph{one} in the sense of~\cite{mr02}) and prove a Gras-type conjecture, relating these Kolyvagin systems to appropriate ideal class groups, refining the results of~\cite{ru92} (in a sense we explain below), and of~\cite{pr-es, r00} applied to our setting.

\end{abstract}

\maketitle

\section*{Introduction}
  B. Howard, B. Mazur and K. Rubin show in~\cite{mr02} that the existence of Kolyvagin systems relies on a cohomological invariant, what they call the Selmer core rank. When the Selmer core rank is one, they determine the structure of the Selmer group completely in terms of a Kolyvagin system. However, when the Selmer core rank is greater than one, not much could be said. In fact, one does not expect a similar result for the structure of the Selmer group in general, as a reflection of the fact that Bloch-Kato conjectures do not in general predict the existence of special elements, but a regulator, to compute the relevant $L$-values.

An example of a core rank greater than one situation arises if one attempts to utilize the Euler system that would come from the Stark elements (whose existence were predicted by K. Rubin~\cite{ru96}). This is what we study in this paper.

Rubin was first to study the Euler system of Stark units in~\cite{ru92}. He proved a Gras-type formula for the $\chi$-isotypic component of a certain ideal class group under certain assumptions on the character $\chi$. These assumptions essentially ensured that the core Selmer rank of $T=\ZZ_p(1)\otimes\chi^{-1}$, in the sense of Definitions 4.8 and 4.1.11 of~\cite{mr02}, is \emph{one}. We prove Theorem A, in the spirit of Gras conjectures and Rubin's prior results (Theorem 4.6 of~\cite{ru92}) for $all$ even, non-trivial characters $\chi$ which are unramified at all primes of $k$ above $p$. The setting in which Rubin~\cite{ru92} places himself is disjoint from ours, however, from the perspective of Kolyvagin system theory Rubin deals with the core Selmer rank \emph{one} situation, whereas in our setting core Selmer rank could be (and in general is) greater then \emph{one}. In that sense, our results may  be regarded as a refinement of the results of~\cite{ru92}.

Before describing our result in detail, we first set some notation. Fix a totally real number field $k$ of degree $r$ over $\QQ$, an algebraic closure $\overline{k}$ of $k$, and a rational prime $p\neq 2$. Let $\chi$ be a character of Gal$(\overline{k}/k)$ into $\ZZ_p^{\times}$ which is trivial on all complex conjugations  inside Gal$(\overline{k}/k)$ (we say that $\chi$ is $even$) and which is unramified at all primes above $p$. We also assume that $\chi$ has finite prime-to-$p$ order, and let $L$ be the fixed field of $\ker(\chi)$ inside $\overline{k}$.  Fix also a finite set $S$ of places of $k$ that does \emph{not} contain $p$, but contains all infinite places $S_{\infty}$, all places $\lambda$ that divide the conductor $f_{\chi}$ of $\chi$. Assume also that $S$ contains at least $r+1$ places, only $r$ of which (namely the infinite places) split completely in $L/k$. We also assume Leopoldt's conjecture (which we often abbreviate as \textbf{LC}) for the number field $L$.

\begin{ithm} [\textbf{A}]
$|A_L^{\chi}| \leq [\wedge^r (O_L^{\times} \otimes \ZZ_p)^{\chi}: \ZZ_p\varepsilon_{L}^{\chi}]$, where $A_L^{\chi}$ is the $p$-part of the $\chi$-isotypic part of the ideal class group of $L$, and $\varepsilon_{L}^{\chi}$ is a Stark element (as in \S 2.1 below) whose existence is predicted in~\cite{ru96}.
\end{ithm}

We note that $\varepsilon_{L}^{\chi}$ is the $\chi$-part of the Stark element of Rubin~\cite{ru96}. Its definition depends on the choice of two disjoint, non-empty finite sets $S$ and $T$ of places of $L$ (see~\cite{ru96} Conjecture $\textup{B}^{\prime}$). Under the assumptions above on the set $S$, it is guaranteed that the order of vanishing of the modified $L$-function $L_{S,T}(s,\chi)$ (\emph{cf.} \cite{ru96} \S1.1 for a definition of the modified $L$-function) at $s=0$ is $r$ and thus, by the defining property of Rubin's Stark element $\varepsilon_{L}^{\chi} \neq 0$. Although the definition of the element $\varepsilon_{L}^{\chi}$ depends on the finite sets $S$ and $T$, the truth of the statement of Theorem~(\textbf{A}) does not as long as $S$ satisfies the hypotheses above. That's why we fix once and for all such $S$ and drop $S$ from the notation wherever it is convenient.

The proof of Theorem (\textbf{A}) has three steps: We first construct in Section~\ref{sec:modified} Selmer structures which we denote by $(T, \FF_{\al})$ modifying the local conditions for the classical Selmer structure $(T, \FF_{cl})$ at $p$. We prove that (Proposition~\ref{modifiedcorerank}) Selmer core rank of each $(T, \FF_{\al})$ is \emph{one}.

This means (by Theorem 5.2.10 of~\cite{mr02}) the $\ZZ_p$-module $\textbf{KS}(T,\FF_{\al})$ of Kolyvagin systems for    $(T, \FF_{\al})$ is free of rank \emph{one}. The second step is to construct these Kolyvagin systems using Stark elements of Rubin~\cite{ru96}.

Once we have these Kolyvagin systems constructed, we apply the Kolyvagin system machinery to obtain bounds for the dual (modified) Selmer groups. In Section~\ref{sec:applications} we compute these bounds quite explicitly in terms of the Stark elements, then compare the modified (dual) Selmer group to the classical (dual) Selmer group and thus obtain the desired bound on the classical (dual) Selmer group. We note that (Proposition~\ref{explicit}) the classical (dual) Selmer group is exactly the relevant part of the ideal class group.

In the very last section, we prove that above inequality is an equality using an analytic class number formula. We still assume \textbf{LC} for $L$. Again the assumptions (and remarks) above for the set $S$ are in effect; and how $T$ should be chosen to guarantee this equality is explained in \S\ref{subsec:stark}.
\begin{ithm}[\textbf{B}]
$|A_L^{\chi}| = [\wedge^r (O_L^{\times})^{\chi}: \ZZ_p\varepsilon_{L}^{\chi}]$
\end{ithm}

One could of course start with Stark elements of Solomon for the values of Twisted-zeta functions at $s=1$ to construct the Kolyvagin systems we need, instead of Rubin's Stark elements. In fact, almost tautologically, one would obtain identical results to what we prove in this paper; by comparing Solomon's Stark elements to Rubin's as in~\cite{solomon-1} Remark 5.3. Possibly more interesting approach would be to start with Solomon's $p$-adic Stark conjectures~\cite{solomon-2}. The author hopes to check in the near future that the ``$p$-adic Stark elements" (i.e. the solutions to $p$-adic Stark conjectures) and the techniques employed in this paper could be used to prove identical results, and also to deduce a relation between ``complex" and ``$p$-adic" Stark conjectures using the rigidity of Kolyvagin systems.

We also remark that Popescu~\cite{popescu} has proved Gras-type conjectures (for arbitrary order of vanishing) for function fields of characteristic $p$ using different techniques then ours.

It's a pleasure to thank Karl Rubin for his guidance, encouragement and  patience. We also thank Christian Popescu and David Solomon for helpful conversations and their valuable remarks. We also thank the anonymous referee for many corrections and suggestions to improve the exposition.

\newpage

\section{Modified Selmer structures for $\ZZ_p(1) \otimes \chi^{-1}$}
\label{sec:modified}
Fix once and for all a totally real number field $k$ of degree $r$ over $\QQ$, an algebraic closure $\overline{k}$ of $k$, and a rational prime $p\neq 2$. Let $\chi$ be an $even$ character of Gal$(\overline{k}/k)$ ($i.e.$ it is trivial on all complex conjugations inside Gal$(\overline{k}/k)$) into $\ZZ_p^{\times}$ that has finite order, and let $L$ be the fixed field of $\hbox{ker}(\chi)$ inside $\overline{k}$.  We define $f_{\chi}$ to be the conductor of $\chi$, and $\Delta$ to be the Galois group $\hbox{Gal}(L/k)$ of the extension $L/k$.

For any abelian group $A$, $A^{\wedge}$ will always denote its $p$-adic completion. Define also $A^{\chi}$ to be the $\chi$-isotypic component of $A^{\wedge}$.

\subsection{Selmer groups for $T$}
Let $T$ be  the Gal$(\overline{k}/k)$-representation $\ZZ_p(1)\otimes~\chi^{-1}$.

One may identify $H^1(L,\ZZ_p(1))$ with $$L^{\times,\wedge} := \lim_{\stackrel{\longleftarrow}{n}}L^{\times}/(L^{\times})^{p^n}$$ using Kummer theory. For any rational prime $\ell$,  define $$H^1(L_{\ell},\ZZ_p(1)):=\bigoplus_{v \mid \ell} H^1(L_v,\ZZ_p(1)) $$  the semi-local Galois cohomology group  at $\ell$, which similarly may be identified by $L_{\ell}^{\times,\wedge}$  where $L_{\ell}=L\otimes\QQ_{\ell}$. Denote by $U_{L,\ell}$ the ($p$-adic completion of the) local units inside  $L_{\ell}^{\times,\wedge}$. Then the classical Selmer structure $\FF_{cl}$ (see Definition 2.1.1 of~\cite{mr02} for the definition of \emph{a Selmer structure}) on $\ZZ_p(1)$ is defined by the local conditions determined by $\{U_{L,\ell} \subset (L\otimes\QQ_{\ell})^{\times,\wedge}\}$.

Likewise,   $$H^1(k,T)=H^1(k,\ZZ_p(1)\otimes{\chi^{-1}})=(H^1(L,\ZZ_p(1))\otimes{\chi^{-1}})^{\Delta}=
(L^{\times})^{\chi}$$ If we
define $H^1(k_{\ell},T)$ analogously, one easily observes that we may identify this semi-local cohomology group with $(L_{\ell}^{\times})^{\chi}$. We similarly define the classical Selmer structure $\FF_{cl}$ on $T$, by the local conditions determined by $\{U_{L,\ell}^{\chi} \subset [(L\otimes\QQ_{\ell})^{\times}]^\chi\}$.

For brevity we denote $U_{L, p}$ by $V_L$. 
Since $\Delta$ has order prime to $p$, $\ZZ_p[\Delta]$ is semi-simple hence it follows that $V_{L}$ is a free $\ZZ_p[\Delta]$-module of rank $r=[k:\QQ]$. Let $\al$ be a free $\ZZ_p[\Delta]$-submodule $V_L$ of rank \emph{one}, such that  the quotient $V_{L}/\al$ is a free $\ZZ_p[\Delta]$-module (of rank $r-1$). Then the $\chi$-~isotypic component $\al^{\chi}$ is free of rank \emph{one} over $\ZZ_p$. Since taking $\chi$-components is exact, it also follows that $V_L^{\chi}/\al^{\chi}$ is free of rank $r-1$ over $\ZZ_p$.

We are now ready to define the modified Selmer structures $(\FF_{\al}, T)$ for $T$, for each choice of $\al$. Let

\begin{itemize}
\item $H^1_{\FF_{\al}}(k_\ell,T):= U_{L,\ell}^{\chi} \subset [(L\otimes\QQ_{\ell})^{\times,\wedge}]^{\chi}=H^1(k_\ell,T)$, if $\ell \neq p$,

\item $H^1_{\FF_{\al}}(k_p,T):= \al^{\chi} \subset U_{L,p}^{\chi} \subset [(L\otimes\QQ_{p})^{\times,\wedge}]^{\chi}= H^1(k_p,T)$.
\end{itemize}

Namely, we shrink the local conditions that defined the classical Selmer structure $(T,\FF_{cl})$ at $p$ down to the $\ZZ_p$-line $\al^{\chi}$ and do nothing for the choices at $\ell\neq p$ to obtain $(T, \FF_{\al})$.

We define the Selmer group for the local conditions  $\FF=\FF_{cl}$ or $\FF=\FF_{\al}$ to be $$H^1_{\FF}(k,T):= \textup{ker}\{H^1(k,T) \lra \prod_{\ell}H^1(k_{\ell},T)/H^1_{\FF}(k_{\ell},T)\}$$ where the product runs over all rational primes $\ell$.

\subsection{Local Duality and Selmer groups for $T^*$}
We now define the dual Selmer structure to $(T,\FF_{\al})$ and the dual Selmer group. Later in this section, we will compute the classical and the modified Selmer groups explicitly and compare their sizes to each other.

Let $$T^*=\textup{Hom}_{\ZZ_p}(T,\QQ_p/\ZZ_p)(1)\cong \QQ_p/\ZZ_p\otimes\chi$$ be the Cartier dual of $T$. For any prime $\lambda$ of $k$, let $<,>_{\lambda}$ denote the local Tate pairing $$<,>_{\lambda}: H^1(k_{\lambda},T) \times H^1(k_{\lambda},T^*) \lra \QQ_p/\ZZ_p$$

\begin{define} For each rational prime $\ell$ define $H^1_{\FF_{\al}^{*}}(k_{\ell},T^*)$ to be the orthogonal complement of $H^1_{\FF_{\al}}(k_{\ell},T)$ under the induced pairing.  The Selmer structure $(T^*,\FF_{\al}^*)$ will be referred to as the dual Selmer structure to $(T,\FF_{\al})$.
\end{define}

Define similarly  $(T^*,\FF_{cl}^*)$, the dual Selmer structure to the classical Selmer structure on $T$.

For $\FF=\FF_{\al}$ or $\FF=\FF_{cl}$, we \emph{propagate} the Selmer structure $(T, \FF)$  (\emph{resp.} $(T^*, \FF^{*})$) to $T/p^nT$ (\emph{resp.} to $T^*[p^n]$) as follows:

Define $H^1_{\FF}(k_{\ell},T/p^nT)$ (\emph{resp.} $H^1_{\FF^*}(k_{\ell},T^*[p^n])$) as the image (\emph{resp.} inverse image) of $H^1_{\FF}(k_{\ell},T)$ (\emph{resp.} $H^1_{\FF^*}(k_{\ell},T^*)$) under the maps induced by $$T \twoheadrightarrow T/p^nT, \, \, \,\,\,T^*[p^n] \hookrightarrow  T^*$$ See also Example 1.1.2 of~\cite{mr02}. We still denote the Selmer structures on $T/p^nT$ (\emph{resp.} on $T^*[p^n]$) obtained this way by $(T/p^nT, \FF)$ (\emph{resp.} by $(T^*[p^n],\FF^*$) for $\FF=\FF_{\al}$ or $\FF=\FF_{cl}$.

\begin{prop}
\label{explicit}
For every positive integer $n$, we have a natural exact sequence  \begin{itemize}
\item $0 \lra (O_L^{\times}/(O_L^{\times})^{p^n})^{\chi} \lra H^1_{\FF_{cl}}(k,T/p^nT)\lra A_L^{\chi}[p^n]\lra 0$

and an isomorphism\\

\item $H^1_{\FF^*_{cl}}(k,T^*[p^n]) \cong \textup{Hom}(A_L^{\chi}, \ZZ/p^n\ZZ)$
\end{itemize}
 where $A_L$ is the ideal class group of L, and $O_L$ is the ring of integers of $L$.

\end{prop}

For a proof, see  \S 2, Proposition 2.6 of \cite{r00}, or \S 6.1 of~\cite{mr02}.

It is easy to see that the classical Selmer structure $(T,\FF_{cl})$ and the modified Selmer structure $(T,\FF_{\al})$ satisfy the hypotheses {H.0-H.5} and (by Lemma 3.7.1 of~\cite{mr02}) {H.6} of \cite{mr02} \S 3.5 (with base field $\QQ$ in their treatment replaced by $k$). Therefore, the existence of Kolyvagin systems for these Selmer structures will be decided by their \emph{core Selmer ranks} (for a definition $cf.$ Definitions 4.1.8 and 4.1.11 of~\cite{mr02}). Let $\XX(T,\FF)$ denote the core Selmer ranks of the Selmer structures $(T,\FF)$ for $\FF=\FF_{cl}$ or for $\FF=\FF_{\al}$. Since the hypotheses {H0-H5} and (by Lemma 3.7.1 of~\cite{mr02}) {H6} hold, $\XX(T,\FF)$ will be (as in Definition 5.2.4 of~\cite{mr02}, using Theorem 4.1.3 of~\cite{mr02}) the common value of $\XX(T/p^nT,\FF)$. Further, it will be quite explicitly given by:
\begin{prop}[compare to Theorem 5.2.5 of~\cite{mr02}]
$$\XX(T,\FF)=\max\{\textup{rank}_{\ZZ_p} (H^1_{\FF}(k,T)) -\textup{corank}_{\ZZ_p}(H^1_{\FF^*}(k,T^*)) ,0\}$$ In particular the Selmer core rank $\XX(T,\FF_{cl})$ is $r=[k:\QQ]$.
\end{prop}
\begin{proof}
The first part is proved exactly in the same way as Theorem~5.2.5 of~\cite{mr02}, which is only stated (and proved) when the base field is $\QQ$. The second part follows now from the first part and Proposition~\ref{explicit}, using the fact that $(O_L^{\times})^{\chi}$ is free of rank $r$ over $\ZZ_p$ (\emph{cf.} Proposition I.3.4 of~\cite{tate}).
\end{proof}

We now describe the modified Selmer groups explicitly. Recall our definition that  $O_L^{\times,\wedge}$ is  the $p$-adically completed group of units  of $L$, and  $(O_L^{\times})^{\chi}$ is the $\chi$-isotypic component of $O_L^{\times,\wedge}$. Recall also that $V_L$ is (the $p$-adic completion of) the local units inside $L\otimes \QQ_p$ and $V_L^{\chi}$ is the $\chi$-isotypic component of $V_L$. We also recall that $V_L$ is a free $\ZZ_p[\Delta]$-module of rank $r$, therefore $V_L^{\chi}$ is a free $\ZZ_p$-module of rank $r$.
\begin{prop}
\label{modifiedexplicit}
Let $\JJ_L$ denote the id\'eles of $L$. Then
\begin{itemize}
\item[(i)]  $H^1_{\FF_{\al}}(k,T) = \textup{ker}\{ (O_L^{\times})^{\chi} \stackrel{\iota}{\lra} V_L^{\chi}/{\al^{\chi}}\} $, where $\iota$ is induced from the natural localization map $(O_L^{\times})^{\chi} \ra V_L^{\chi}$ (which will also be denoted by $\iota$ later),\\
\item[(ii)] $H^1_{\FF_{\al}^{*}}(k,T^*) \cong \textup{Hom}(\JJ_L/\overline{L^{\times}(\prod_{v\mid\infty}L_v^{\times}\prod_{v\nmid p} U_v  {\al})},\QQ_p/\ZZ_p)^{\chi^{-1}},$\\ where $U_v$ denotes the local units inside the completion $L_v$ of $L$.
\end{itemize}
\end{prop}

\begin{proof}
By the definition of $\FF_{\al}$ and by Proposition~\ref{explicit}, there is a commutative diagram
$$\xymatrix{
0 \ar[r] & H^1_{\FF_{\al}}(k,T) \ar[r]\ar@{=}[d]& H^1_{\FF_{cl}}(k,T) \ar[r]^{\iota}\ar@{=}[d]&\frac{ H^1_{\FF_{cl}}(k_p,T)}{H^1_{\FF_{\al}}(k_p,T)} \ar@{=}[d]\\
0\ar[r] & H^1_{\FF_{\al}}(k,T) \ar[r]&(O_L^{\times})^{\chi} \ar[r]^{\iota}&V_L^{\chi}/\al^{\chi}
}$$ with exact rows. This proves (i).

(ii) also follows easily noting that $H^1(k,T^*)=\textup{Hom}(G_L,\QQ_p/\ZZ_p)^{\chi^{-1}}$ and by class field theory.
\end{proof}

We will now compare the classical (dual) Selmer group (which turned out to be the $\chi$-part of the ideal class group of $L$) to our modified (dual) Selmer groups.

\begin{prop}
\label{comparison}
$H^1_{\FF_{cl}^*}(k,T^*) \subset H^1_{\FF_{\al}^*}(k,T^*)$ and $$[H^1_{\FF_{\al}^*}(k,T^*):H^1_{\FF_{cl}^*}(k,T^*)]=[V_L^{\chi} :\al^{\chi}\cdot \iota((O_L^{\times})^{\chi})].$$ \end{prop}

\begin{proof}
By definitions of the Selmer groups $H^1_{\FF_{cl}} (k, T)$ and $H^1_{\FF_{\al}} (k, T)$ (and the dual Selmer groups) there are exact sequences
$$\xymatrix{
0 \ar[r] &H^1_{\FF_{\al}} (k, T) \ar[r]& H^1_{\FF_{cl}} (k, T) \ar[r]^(.40) {\iota}& \frac{H^1_{\FF_{cl}} (k_p, T)}{H^1_{\FF_{\al}} (k_p, T)}=V_L^{\chi}/\al^{\chi}\\
0 \ar[r]& H^1_{\FF_{cl}^{*}} (k, T^*) \ar[r]& H^1_{\FF_{\al}^{*}} (k, T^*) \ar[r]^{\iota^{*}} &\frac{H^1_{\FF_{\al}^*}(k_p,T^{*})}{H^1_{\FF_{cl}^*}(k_p,T^{*})}
}$$
Further, by Poitou-Tate global duality, the image of $\iota$ and $\iota^{*}$ are orthogonal complements under the local Tate pairing  ($cf.$ Theorem 2.3.4 of~\cite{mr02} or Theorem 1.7.3(i) of~\cite{r00}). It follows then that (see Theorem 1.7.3(ii) of~\cite{r00}) $$\xymatrix{H^1_{\FF_{\al}^{*}} (k, T^*)/H^1_{\FF_{cl}^{*}} (k, T^*)\ar[r]^{\sim} &\textup{Hom}(\textup{coker}(\iota),\QQ_p/\ZZ_p)}$$
By Proposition~\ref{explicit} it follows immediately that  $$\textup{coker}(\iota) \cong V_L^{\chi}/\al^{\chi}\iota((O_L^{\times})^{\chi}).$$ 
This proves that $$[H^1_{\FF_{\al}^*}(k,T^*):H^1_{\FF_{cl}^*}(k,T^*)]=|\textup{coker}(\iota)|=[V_L^{\chi} :\al^{\chi}\cdot \iota((O_L^{\times})^{\chi})].$$
\end{proof}

\begin{rem}
We note that the equality Proposition~\ref{comparison} still holds true even if one of the indices is infinite, in the sense that if one side is infinite, the other is too.
\end{rem}
\begin{rem}
\label{rem:LCshows}
Assume in this paragraph that Leopoldt's conjecture holds for $L$. Note that $$V_L^{\chi} /\al^{\chi} \iota((O_L^{\times})^{\chi})\cong V_L^{\chi}\big{/}\iota((O_L^{\times})^{\chi})\Big{/}\al^{\chi}\big{/}\iota((O_L^{\times})^{\chi})\cap\al^{\chi}$$
and thus, by Proposition~\ref{comparison} $$[H^1_{\FF_{\al}^*}(k,T^*):H^1_{\FF_{cl}^*}(k,T^*)]=\frac{[V_L^{\chi}:\iota((O_L^{\times})^{\chi})]}{[\al^{\chi}:\iota((O_L^{\times})^{\chi})\cap\al^{\chi}]}.$$ Note that both the numerator and the denominator are finite since we assumed Leopoldt's conjecture. 
\end{rem}

\begin{prop}
\label{modifiedcorerank}
The Selmer structure $(T,\FF_{\al})$ has Selmer core rank $one$ (see 4.1.8 and 4.1.11 of~\cite{mr02} for a definition of Selmer core rank).
\end{prop}
\begin{proof}
By Proposition 1.6 of ~\cite{wiles}
\begin{align*}
\textup{length}(H^1_{\FF_{\al}}(k,T/p^nT))&-\textup{length}(H^1_{\FF_{\al}^*}(k,T^*[p^n]))=\\
=\textup{length}(H^0(k,T/p^nT))&-\textup{length}(H^0(k,T^*[p^n]))\\
-\sum_{\ell \mid f_{\chi}p} \{\textup{length}(H^0(k_{\ell},T/p^nT)&-\textup{length}(H^1_{\FF_{\al}}(k_{\ell},T/p^nT))\}\\
\end{align*}
which is $p^n\XX(T,\FF_{\al})$. Applying the same formula to $(T,\FF_{cl})$ we see that \begin{align*}
p^n(\XX(T,\FF_{cl})-&\XX(T,\FF_{\al}))=\\
&\textup{length}(H^1_{\FF_{cl}}(k_{p},T/p^nT))-\textup{length}(H^1_{\FF_{\al}}(k_{p},T/p^nT))
\end{align*}
 and this equals $(r-1)p^n$. Since we already know that $\XX(T,\FF_{cl})=r$ the Proposition follows.

\end{proof}
Let  $\textup{\textbf{KS}}(T,\FF_{\al})$ denote the $\ZZ_p$-module of Kolyvagin systems for the Selmer structure $(T,\FF_{\al})$. See Definition 3.1.3 of~\cite{mr02} for a precise definition.
\begin{cor}
 $\textup{\textbf{KS}}(T,\FF_{\al})$ is free of rank \emph{one} as a $\ZZ_p$-module, generated by a Kolyvagin system $\kappa \in \textup{\textbf{KS}}(T,\FF_{\al})$ whose image inside of $\textup{\textbf{KS}}(T/pT,\FF_{\al})$ is nonzero.
\end{cor}

\begin{proof}
This is immediate after Proposition~\ref{modifiedcorerank} and Theorem 5.2.10 of~\cite{mr02}.
\end{proof}

\newpage
\section{Kolyvagin systems of Stark units}
\label{sec:kolsys}

In this section we review Rubin's integral refinement of Stark conjectures and construct Kolyvagin systems for the modified Selmer structures $(T,\FF_{\al})$ using Stark elements of Rubin.

For the rest of the paper we assume this refined version of the Stark conjectures (Conjecture $\textup{B}^\prime$ of~\cite{ru96}).

We first set some notation. Assume $k,\chi, f_{\chi}$ and $L$ are as above. For a cycle $\tau$ of the number field $k$ let $k(\tau)$ be the maximal $p$-extension inside the ray class field of $k$ modulo $\tau$. Define $L(\tau)$ to be the composite of $k(\tau)$ and $L$. Let $$\mathcal{K} =\{L(\tau): \tau \hbox{ is a (finite) cycle of } k \hbox{ prime to } f_{\chi}p\}$$ be a collection of abelian extensions  of $k$,  where $k(\tau)$, $L$, $L(\tau)$ and $f_{\chi}$ are defined above.

\subsection{Stark units of Rubin and Euler systems}
\label{subsec:stark}
Fix a finite set $S$ of places of $k$ that does \emph{not} contain any prime above $p$, but contains all infinite places $S_{\infty}$ and all places $\lambda$ that divide the conductor $f_{\chi}$ of $\chi$. Assume that $|S| \geq r+1$. For each $K \in \kk$ let $S_K=S \cup \{\hbox{places of k at which K is ramified}\}$ be another set of places of $k$. Let $O_{K,S_K}^{\times}$ denote the $S_K$ units of $K$, and $\Delta_K$ (\emph{resp.} $\delta_K$) denote $\hbox{Gal}(K/k)$ (\emph{resp.} $|\hbox{Gal}(K/k)|$).  Conjecture $\textup{B}^{\prime}$ of~\cite{ru96} predicts the existence of certain elements\footnote{In fact, Rubin's conjecture predicts that these elements should be inside  $\frac{1}{\delta_K}{\wedge^r} O_{K,S_K,T}^{\times}$ where $T$ is a finite set of primes disjoint from $S_K$ chosen so that the group $O_{K,S_K,T}^{\times}$ of units which are congruent to 1 modulo all the primes in $T$ is torsion-free. However, in our case any set $T$ which contains a prime other than 2 will suffice (since all the fields that appear in our paper are totally real). Further, $T$ may be chosen in a way that (e.g. $T=\{p\}$) the extra factors that appear in the definition of the modified $zeta\,function$ for $K$ ($c.f.$ \S1 of~\cite{ru96} for more detail on these zeta functions) will be prime to $p$, when they are evaluated at $0$. We note that for such a chosen $T$, we have $O_{K,S_K,T}^{\times,\wedge}=O_{K,S_K}^{\times,\wedge}$ , for example by the exact sequence (1) in~\cite{ru96}. Since in our paper we will work with the $p$-adic completion of the units, we will safely exclude $T$ from our notation.} $$\varepsilon_{K,S_K} \in \Lambda_{K,S_K} \subset \frac{1}{\delta_K}{\wedge^r} O_{K,S_K}^{\times}$$ where $\Lambda_{K,S_K}$ is defined in \S 2.1 of~\cite{ru96} and has the property that for any homomorphism $$\tilde{\psi} \in \hbox{Hom}_{\QQ_p[\Delta_K]}(\wedge^{r}O_{K,S_K}^{\times,\wedge} \otimes \QQ_p,O_{K,S_K}^{\times, \wedge} \otimes \QQ_p)$$ that is induced from a homomorphism $$\psi \in \hbox{Hom}_{\ZZ_p[\Delta_K]}(\wedge^{r}O_{K,S_K}^{\times,\wedge},O_{K,S_K}^{\times, \wedge})$$ we have $\tilde{\psi}(\Lambda_{K,S_K}) \subset O_{K,S_K}^{\times, \wedge}$. We note that the $r$-th exterior power ${\wedge^r} O_{K,S_K}^{\times}$ (and other exterior powers which appear below) is taken in the category of $\ZZ_p[\Delta_K]$-modules. The elements $\varepsilon_{K,S_K}$ (which we call Stark elements) satisfy the distribution relation to be satisfied by an Euler system (Proposition 6.1 of~\cite{ru96}). We denote the image of $\varepsilon_{K,S_K}$ inside the $\ZZ_p$-module $\frac{1}{\delta_K}\wedge^r O_{K,S_K}^{\times,\wedge}$ also by $\varepsilon_{K,S_K}$. Since $S$ is fixed (therefore $S_K$, too), we will often drop $S$ or $S_K$ from notation and denote $\varepsilon_{K,S_K}$ by $\varepsilon_{K}$, or sometimes use $S$ instead of $S_K$ and denote $O_{K,S_K}$ by $O_{K,S}$.

As before, let $V_K$ be the $p$-adic completion of the units  in $K\otimes \QQ_p$. Any $\phi_K \in \wedge^{r-1}\hbox{Hom}_{\ZZ_p[\Delta_K]}(V_K, \ZZ_p[\Delta_K])$ induces a homomorphism, which we still denote by $\phi_K$, inside $\hbox{Hom}_{\ZZ_p[\Delta_K]}(\wedge^{r}O_{K,S}^{\times,\wedge},O_{K,S}^{\times, \wedge})$, as discussed in \S 1.2 and \S 6.3 of~\cite{ru96}. Let  $M$ be be the composite of the fields inside $\kk$, let $$\ZZ_p[[\textup{Gal}(M/k)]]:=\lim_{\stackrel{\longleftarrow}{K \in \kk}} \ZZ_p [\textup{Gal}(K/k)]$$ denote the completed group ring of \textup{Gal}(M/k) and $$\vv:=\lim_{\stackrel{\longleftarrow}{K\in\kk}}V_K$$ We define $$\bigwedge^{r-1}\textbf{Hom}(\mathbb{V},\ZZ_p[[\textup{Gal}(M/k)]]):=\lim_{\stackrel{\longleftarrow}{K \in \kk}} \bigwedge^{r-1}\hbox{Hom}_{\ZZ_p[\Delta_K]}(V_K, \ZZ_p[\Delta_K]) $$ where the inverse limit is with respect to the natural maps induced from the inclusion map $V_K \hookrightarrow V_{K^{\prime}}^{\hbox{\tiny Gal}(K^{\prime}/K)}$ and the isomorphism
 \begin{align*}
  \ZZ_p[\Delta_{K^{\prime}}]^{\hbox{\tiny Gal}(K^{\prime}/K)} & \tilde{\lra}\ZZ_p[\Delta_K]\\
\mathbf{N}^{K^{\prime}}_K &\longmapsto 1
 \end{align*} for $K \subset K^{\prime}$.

 \begin{prop}
 \label{krasner-2}
 For any $K \in \kk$ the projection map $$\hhh \lra \bigwedge^{r-1}\textup{Hom}_{\ZZ_p[\Delta_K]}(V_K, \ZZ_p[\Delta_K])$$ is surjective.

\end{prop}

\begin{proof}
Let $K \in \kk$. The quotient of the local units by the local 1-units has order prime
to $p$, because it injects (under reduction) into the multiplicative
group of the residue field(s).  So the $p$-adic completion of the local
units (which we have denoted by $V_K$) is the same as the $p$-adic completion of the local 1-units. The rest of the proof follows at once from Corollary 6.5 of~\cite{ru96}. 
\end{proof}
 Since $S_K$ contains no primes above $p$, $O_{K,S_K}^{\times,\wedge}$ maps canonically to $V_K$ via a $\ZZ_p[\Delta_K]$-equivariant map. This induces a natural map \begin{equation}\label{eqn:view}\lim_{\stackrel{\longleftarrow}{K \in \kk}} \bigwedge^{r-1}\hbox{Hom}_{\ZZ_p[\Delta_K]}(V_K, \ZZ_p[\Delta_K]) \lra \lim_{\stackrel{\longleftarrow}{K \in \kk}} \bigwedge^{r-1}\hbox{Hom}_{\ZZ_p[\Delta_K]}(O_{K,S_K}^{\times,\wedge}, \ZZ_p[\Delta_K])\end{equation} The image of $\Phi \in \hhh$ under the above map will still be denoted by $\Phi$, when there is no confusion.

\begin{prop}
\label{eulersystem}
Let $\{\phi_K\}=\Phi \in \hhh$ and let $$\varepsilon_{K,S_K,\Phi}=\phi_K(\varepsilon_{ K,S_K}) \in O_{K,S_K}^{\times,\wedge}$$ \textup{(}Here we view $\phi_K$ as an element of\,  $\textup{Hom}_{\ZZ_p[\Delta_K]}(\wedge^rO_{K,S_K}^{\times,\wedge}, O_{K,S_K}^{\times,\wedge})$ via \textup{(}\ref{eqn:view}\textup{)} above and the map \textup{(}4\textup{)} (with $k=r-1$) of~\cite{ru96}.\textup{)} Then the collection $\{\varepsilon_{K,S_K,\Phi}\}_{K \in \kk}$ is an Euler system for the $\textup{Gal}(\overline{k}/k)$-representation $\ZZ_p(1)$ in the sense of Definition 2.1.1 of~\cite{r00} \textup{(}with condition \textup{(ii)} replaced by $\textup{(ii)}^{\prime} \textup{(b)}$ in \S 9.1 of~\cite{r00}\textup{)}.
\end{prop}

\begin{proof}
This is Proposition 6.6 of~\cite{ru96}. We only remark that since $\varepsilon_{K,S_K} \in \Lambda_{K,S_K}$ it follows that $\varepsilon_{K,S_K,\Phi} \in O_{K,S_K}^{\times,\wedge} \subset K^{\times,\wedge}$ and Kummer theory identifies $K^{\times,\wedge}$  with $H^1(K,\ZZ_p(1))$.
\end{proof}

\subsection{Twisting by the character $\chi$}
Let $\chi$ be a character as defined in Section~\ref{sec:modified}. In this section we will show how to twist the Euler systems we obtain the previous section by the character $\chi$. However, one should notice that in Section~\ref{sec:applications} we will be using the Euler system for the twist by $\chi^{-1}$.

Let $\eee_{\chi}$ denote the idempotent $\frac{1}{|\Delta|}\sum_{\delta \in \Delta}\chi(\delta)\delta^{-1}$. For any finite cycle $\tau$ which is prime to $pf_{\chi}$ we define
\begin{equation}
\label{twist eqn}
\begin{array}{rcl}
\varepsilon_{L(\tau)}^{\chi}=\epsilon_{\chi}\varepsilon_{L(\tau),S} &\in&  \frac{1}{|G_{\tau}|} \epsilon_{\chi}\wedge^r O_{L(\tau),S}^{\times,\wedge}\\\\
&=&\frac{1}{|G_{\tau}|} \wedge^r (O_{L(\tau),S}^{\times})^{\chi}
\end{array}
\end{equation}                                                                                                   where $G_{\tau}:=\hbox{Gal}(k(\tau)/k)$, which is the $p$-part of $\Delta_{\tau}:=\hbox{Gal}(L(\tau)/k)\cong G_{\tau} \times \Delta$ (since $|\Delta|$ is prime to $p$, it follows that $L$ and $k(\tau)$ are linearly disjoint over $k$). We also note that the equality on the second line above \begin{equation}\label{wedge eqn}(\wedge^r O_{L(\tau),S}^{\times})^{\chi}=\epsilon_{\chi}\wedge^r O_{L(\tau),S}^{\times,\wedge} =\wedge^r \epsilon_{\chi}O_{L(\tau),S}^{\times,\wedge}=\wedge^r(O_{L(\tau),S}^{\times})^{\chi}\end{equation} holds simply because $\epsilon^{r}_{\chi}=\epsilon_{\chi}$ (we recall once again that the exterior products in (\ref{twist eqn}) and (\ref{wedge eqn}) are taken in the category of $\ZZ_p[\Delta_{\tau}]$-modules).
\begin{lemma}
\label{twisted-arguments}
For any $\{\phi_K\}=\Phi \in \hhh$ we have $$\varepsilon_{L(\tau),\Phi}^{\chi}:=\phi_{L(\tau)}(\varepsilon_{L(\tau)}^{\chi}) \in (O_{L(\tau),S}^{\times})^{\chi}$$
\end{lemma}

\begin{proof}
$\varepsilon_{L(\tau)}^{\chi}=\epsilon_{\chi}\varepsilon_{L(\tau),S}$, and by definition, $\phi_{L(\tau)}$ is $\ZZ_p[\Delta_{\tau}]$-equivariant, therefore $$\phi_{L(\tau)}(\varepsilon_{L(\tau)}^{\chi})=\epsilon_{\chi}\phi_{L(\tau)}(\varepsilon_{L(\tau),S}) \in \wedge^r(O_{L(\tau),S}^{\times})^{\chi}$$
\end{proof}

The natural inclusion $$V_{L(\tau)}^{\chi} \hookrightarrow V_{L(\tau)}$$ induces a map
$$
\bigwedge^{r-1}\hbox{Hom}_{\ZZ_p[\Delta_{\tau}]}(V_{L(\tau)}, \ZZ_p[\Delta_{\tau}]) \lra  \bigwedge^{r-1}\hbox{Hom}_{\ZZ_p[G_{\tau}]}(V_{L(\tau)}^{\chi}, \ZZ_p[\Delta_{\tau}]^{\chi})   $$
\begin{lemma}
\label{homs}
The map $$\Gamma: \bigwedge^{r-1}\textup{Hom}_{\ZZ_p[\Delta_{\tau}]}(V_{L(\tau)}, \ZZ_p[\Delta_{\tau}]) \lra  \bigwedge^{r-1}\textup{Hom}_{\ZZ_p[G_{\tau}]}(V_{L(\tau)}^{\chi}, \ZZ_p[\Delta_{\tau}]^{\chi})$$
is surjective.
\end{lemma}

\begin{proof}
Using the decompositions $$V_L={\bigoplus}_{\psi \in \Hat{\Delta}} V_L^{\psi}=V_L^{\chi}\oplus (V_L^{\chi})^{\perp}$$ and $$\ZZ_p[\Delta_{\tau}]={\bigoplus}_{\psi \in \Hat{\Delta}} \ZZ_p[\Delta_{\tau}]^{\psi}=\ZZ_p[\Delta_{\tau}]^{\chi}\oplus (\ZZ_p[\Delta_{\tau}]^{\chi})^{\perp}$$ into $\ZZ_p[\Delta_{\tau}]$-modules  the proposition follows at once.
\end{proof}

\begin{rem}
\label{free}
For each $\tau$, fixing a basis $\{\eee_{\chi}^{L(\tau)}\}$, where $$\eee_{\chi}^{L(\tau)}=\frac{1}{|\Delta|}\sum_{\delta \in \Delta} \chi^{-1}(\delta)\delta \in \ZZ_p[\Delta_{\tau}]$$ (we view $\eee_{\chi}^{L(\tau)}$ as an element of $\ZZ_p[\Delta_{\tau}]$ using the decomposition $\Delta_{\tau}=G_{\tau}\times \Delta$) for the free-of-rank-one $\ZZ_p[G_{\tau}]$-module $\ZZ_p[\Delta_{\tau}]^{\chi}$, one may identify $\ZZ_p[\Delta_{\tau}]^{\chi}$ by  ${\ZZ_p[G_{\tau}]}$. Abusing this fact, we will allow ourselves to alternate between  $\ZZ_p[\Delta_{\tau}]^{\chi}$ and  ${\ZZ_p[G_{\tau}]}$, mostly in favor of the latter.
\end{rem}
 Let $$\hhhk:=\lim_{\stackrel{\longleftarrow}{L(\tau) \in \kk}} \bigwedge^{r-1}\textup{Hom}_{\ZZ_p[G_{\tau}]}(V_{L(\tau)}^{\chi}, \ZZ_p[\Delta_{\tau}]^{\chi})$$ The inverse limit is with respect to the natural maps above which we used to define $\hhh$.
\begin{rem}
\label{rem:surjective}
Lemma~\ref{homs} shows that the map $$\hhh \lra \hhhk$$ is surjective.
\end{rem}

\begin{prop}
\label{liftinghoms}
For each $\tau_0$, the projection $$\hhhk \lra  \bigwedge^{r-1}\textup{Hom}_{\ZZ_p[G_{\tau_{0}}]}(V_{L(\tau_{0})}^{\chi}, \ZZ_p[\Delta_{\tau_0}]^{\chi})$$ is surjective.
\end{prop}

\begin{proof}
Proof of this is identical to the proof of Proposition~\ref{krasner-2}.

\end{proof}
\begin{prop}
\label{twisted-euler-system}
Let $\{\phi_{\tau}\}=\Phi \in \hhhk$, then the collection $\{\phi_{\tau}(\varepsilon_{L(\tau)}^{\chi})\}$ is an Euler system for $(T, f_{\chi}p)$ in the sense of Definition 2.1.1 of~\cite{r00} (with condition \textup{(ii)} replaced by $\textup{(ii)}^{\prime} \textup{(b)}$ in \S 9.1 of~\cite{r00}).
\end{prop}
\begin{rem}
We will often write $\varepsilon_{L(\tau),\Phi}^{\chi}$ for $\phi_{\tau}(\varepsilon_{L(\tau)}^{\chi})$. Note that $$\varepsilon_{L(\tau),\Phi}^{\chi} \in (O_{L(\tau),S}^{\times})^{\chi}$$ and thus may also be viewed as an element of $H^1(k(\tau),T)$.
\end{rem}

\begin{proof}
Let $\tilde\Phi$ denote an element of $\hhh$ that lifts $\Phi$. Such a $\tilde\Phi$ exists by Remark~\ref{rem:surjective}. Then  we have the following commutative diagram:

$$\xymatrix  @C=.28in{
\varepsilon_{L(\tau),S} \ar @{|->}[dd] \ar @/^2pc/@{|->}[rrrrrr]&**[l]\in&\bigwedge^r \Lambda_{L(\tau),S} \ar[rr]^{\tilde{\phi}_{\tau}} \ar[dd]^{\epsilon_{\chi}}&&  **[r]O_{L(\tau),S}^{\times,\wedge} \ar[dd]^{\epsilon_{\chi}} &**[l] \ni&\varepsilon_{L(\tau),S, \tilde{\Phi}} \ar @{|->}[dd]\\
&&&&&\\
\varepsilon_{L(\tau)}^{\chi}\ar @/_2pc/@{|->}[rrrrrr]&**[l]\in&\bigwedge^r \Lambda_{L(\tau),S}^{\chi} \ar[rr]^{\phi_{\tau}} && **[r](O_{L(\tau),S}^{\times})^{\chi}&**[l] \ni &\varepsilon_{L(\tau),\Phi}^{\chi}
}$$
\\\\
What this diagram essentially says is $$\phi_{\tau}(\varepsilon_{L(\tau)}^{\chi})=\epsilon_{\chi}\varepsilon_{L(\tau),S,\tilde{\Phi}} (:=\varepsilon_{L(\tau),\tilde{\Phi}}^{\chi})$$ for any lift $\tilde{\Phi}$ of $\Phi$.

By Proposition~\ref{eulersystem} $\{\varepsilon_{F,S, \tilde{\Phi}}\}_{_{ F\in\kk}}$ is an Euler system for $\ZZ_p(1)$. This means, by Proposition 4.2 and Lemma 4.3 in \S 2 of~\cite{r00}, that $\{\epsilon_{\chi}\varepsilon_{L(\tau),S,\tilde{\Phi}} =\varepsilon_{L(\tau),\tilde{\Phi}}^{\chi}\}_{_{\tau}}$ is an Euler system for $\ZZ_p(1)\otimes{\chi^{-1}}$. But $\varepsilon_{L(\tau),\tilde{\Phi}}^{\chi}=\phi_{\tau}(\varepsilon_{L(\tau)}^{\chi})=\varepsilon_{L(\tau),\Phi}^{\chi}$ and the Proposition follows.
\end{proof}

\begin{define}
\label{def:eulersystem}
Let $\varepsilon_{L(\tau),\Phi}^{\chi}$ be as above. We call the collection $\{\varepsilon_{L(\tau),\Phi}^{\chi}\}_{_{\tau}}$ the Euler system of $\Phi$-Stark elements for $T=\ZZ_p(1)\otimes\chi^{-1}$.
\end{define}

\begin{rem}
\label{integrality}
By definition $\varepsilon_{L}^{\chi}$ belongs to $\frac{1}{|\Delta|} \wedge^{r}(O_{L,S}^{\times})^{\chi}$, where $\Delta=\textup{Gal}(L/k)$. Since $|\Delta|$ is prime to $p$, we have $$\frac{1}{|\Delta|} \wedge^{r}(O_{L,S}^{\times})^{\chi}= \wedge^{r}(O_{L,S}^{\times})^{\chi}$$ hence $$\varepsilon_{L}^{\chi} \in \wedge^{r}(O_{L,S}^{\times})^{\chi}$$ Further let $$\tilde{U}_{L,S}=\{u \in O_{L,S}^{\times,\wedge}: \eee_{\chi_0}u=0 \hbox{ for the trivial character } \chi_0 \in \Hat{\Delta} \}$$ be as in \S 6.3 of~\cite{ru96}. Note that $(O_{L,S}^{\times})^{\chi} \subset \tilde{U}_{L,S}$ since $\chi$ is non-trivial. Also, the proof of Proposition 6.2 of~\cite{ru96} shows that $$\tilde{U}_{L,S} \subset O_{L}^{\times,\wedge}: p\hbox{-adic completion of the global units,}$$ therefore $$\varepsilon_{L}^{\chi} \in  \wedge^{r}(O_{L}^{\times})^{\chi}$$
\end{rem}

\subsection{Choosing the Homomorphisms $\Phi$:}\label{hom}
Recall that $M$ was defined to be the composite of all the fields that are inside the collection $\kk$ and $$\ZZ_p[[\textup{Gal}(M/k)]]:=\lim_{\stackrel{\longleftarrow}{K\in \kk}} \ZZ_p[\textup{Gal}(K/k)]$$ as the completed group ring of $\textup{Gal}(M/k)$. Recall also that  $M_0$ is the fixed field $M^{\Delta}$ of $\Delta$ inside $M$.

Until the end of the paper, we assume that $L/k$ is unramified at all primes of $k$ above $p$. Note then that for any $K \in \kk$, $K/k$ is unramified at all primes of $k$ which are above $p$, and are all totally real. Therefore Krasner's lemma~\cite{krasner} on the 1-units implies:

\begin{lemma}[Krasner]
\label{krasner}
$U_{K, p}$ is a free $\ZZ_p[\textup{Gal}(K/k)]$-module of rank $r=[k:\QQ]$, where $U_{K, p}=\bigoplus_{\wp \mid p} U_{K_{\wp}}$ is the $p$-adic completion of the group of units in $K\otimes\QQ_p$.
\end{lemma}

As an immediate consequence of Krasner's lemma (Lemma~\ref{krasner}) we have
\begin{cor}
\label{thick-line}
$$\vv=\lim_{\stackrel{\longleftarrow}{K\in\kk}}V_K$$ is free of rank $r$ over $\ZZ_p[[\textup{Gal}(M/k)]]$.
\end{cor}

Using Corollary~\ref{thick-line} we may choose a $\ZZ_p[[\textup{Gal}(M/k)]]$-line $\mathbb{L}$  inside $\vv$ such that the quotient $\vv/\mathbb{L}$ is also a free $\ZZ_p[[\textup{Gal}(M/k)]]$-module (of rank $r-1$).
\begin{define}
\label{lines}
For all $L(\tau)=K \in \kk$ let $\al_{K}$ be the image of $\mathbb{L}$ under the (surjective) projection map \, $\vv \ra V_K$. When $\tau=1$ (i.e. when $K=L$), denote $\al_K$ by only $\al$.
\end{define}

Note that $\al_K$ are free $\ZZ_p[\textup{Gal}(K/k)]$-modules of rank \emph{one} for all $K \in \kk$, and that $(\al_K)^{\textup{Gal}(K/K^{\prime})}=\al_{K^{\prime}}$ for all $K^{\prime} \subset K$.

We will often denote $\al_{L(\tau)}$ by simply $\al_{\tau}$.
\begin{rem}
\label{splits}
If there is a prime $\wp$ of $k$ above  $p$ of degree one in $k/\QQ$, then for each $K \in \kk$ there is a natural choice for a $\ZZ_p[\textup{Gal}(K/k)]$-line $\al_K$: Define $$\al_K=\bigoplus_{\hbox{\small{\frakfamily q}}\mid\wp}U_{K_{\hbox{{\small\frakfamily q}}}}$$ where the direct sum is over the places {\frakfamily q} of $K$ that are above $\wp$, and $U_{K_{\hbox{{\small\frakfamily q}}}}$ is defined as before. For $K^{\prime} \subset K$  it is evident that $(\al_{K})^{\textup{Gal}(K/K^{\prime})}=\al_{K^{\prime}}$. We also know by (a variant of) Lemma~\ref{krasner} that $\al_{K}$ is a free $\ZZ_p[\hbox{Gal}(K/k)]$-module of rank $one$. It also follows that the $\chi$-isotypic component $(\al_{L(\tau)})^{\chi}$ of $\al_{L(\tau)}$ is free of rank $one$ over $\ZZ_p[G_{\tau}]$, for $K=L(\tau)$.
\end{rem}

\begin{define}
\label{hli}
We say that an element $$\{\phi_{\tau}\}_{_{\tau}}=\Phi \in \hhhk$$ satisfies $\hli$ if for any $L(\tau)\in \kk$ one has $\phi_{\tau}(\wedge^{r} V_{L(\tau)}^{\chi}) \subset \al_{\tau}^{\chi}$.
\end{define}

We will next construct a specific element $$\Phi_0 \in \hhhk$$ that satisfies $\hli$.

In what follows, we will identify the free of rank $one$ $\ZZ_p[[\textup{Gal}(M_0/k)]]$-module $\ZZ_p[[\textup{Gal}(M/k)]]^{\chi}$ with  $\ZZ_p[[\textup{Gal}(M_0/k)]]$ (as in Remark~\ref{free}), and we allow ourselves to alternate between these two notations.

Fix a basis $$\{\Psi_{\mathbb{L}}^{(i)}\}_{i=1, \dots,r-1}$$ of the free (of rank $r-1$) $\ZZ_p[[\textup{Gal}(M_0/k)]]$-module $$\textup{Hom}_{\ZZ_p[[\textup{Gal}(M_0/k)]]}(\vv^{\chi}/\all^{\chi},\ZZ_p[[\textup{Gal}(M_0/k)]])$$ This then fixes  a basis $\{\psi_{\al_{\tau}}^{(i)}\}_{i=1}^{r-1}$ for the free (of rank r-1) $\ZZ_p[G_{\tau}]$-module $\textup{Hom}_{\ZZ_p[G_{\tau}]}(V_{L(\tau)}^{\chi}/\al_{\tau}^{\chi}, \ZZ_p[G_{\tau}])$ for all $L(\tau) \in \kk$; such that $\{\psi_{\al_{\tau}}^{(i)}\}_{\tau}$ are compatible with respect to the surjections $$\textup{Hom}_{\ZZ_p[G_{\tau}]}(V_{L(\tau)}^{\chi}/\al_{\tau}^{\chi},\ZZ_p[G_{\tau}]) \lra \textup{Hom}_{\ZZ_p[G_{\tau^{\prime}}]}(V_{L(\tau^{\prime})}^{\chi}/\al_{\tau^{\prime}}^{\chi},\ZZ_p[G_{\tau^{\prime}}])$$ for all $\tau^{\prime} \mid \tau$. Note that the homomorphism $$\bigoplus_{i=1}^{r-1}\psi_{\al_{\tau}}^{(i)}:V_{L(\tau)}^{\chi}/\al_{\tau}^{\chi} \lra \ZZ_p[G_{\tau}]^{r-1}$$ is an isomorphism of $\ZZ_p[G_{\tau}]$-modules, for all $\tau$.

Let $\psi_{\tau}^{(i)}$ denote the image of $\psi_{\al_{\tau}}^{(i)}$ under the canonical injection $$\textup{Hom}_{\ZZ_p[G_{\tau}]}(V_{L(\tau)}^{\chi}/\al_{\tau}^{\chi},\ZZ_p[G_{\tau}]) \hookrightarrow \textup{Hom}_{\ZZ_p[G_{\tau}]}(V_{L(\tau)}^{\chi},\ZZ_p[G_{\tau}])$$ Note then that $$\Psi_{\tau}:=\bigoplus^{r-1}_{i=1}\psi_{\tau}^{(i)}: V_{L(\tau)}^{\chi} \lra \ZZ_p[G_\tau]$$ is surjective and $\textup{ker}(\Psi_{\tau})=\al_{\tau}^{\chi}$.


Define $$\phi_{\tau} = \psi_{\tau}^{(1)}\wedge\psi_{\tau}^{(2)} \wedge\dots \wedge\psi_{\tau}^{(r-1)} \in \bigwedge^{r-1}\textup{Hom}_{\ZZ_p[\Delta_{\tau}]}(V_{L(\tau)}^{\chi},\ZZ_p[G_{\tau}]) $$  Note once again that for $\tau^{\prime}|\tau$, $\phi_{\tau}$ maps to $\phi_{\tau^{\prime}}$ under the surjective (by Proposition~\ref{liftinghoms}) homomorphism $$\bigwedge^{r-1}\textup{Hom}_{\ZZ_p[\Delta_{\tau}]}(V_{L(\tau)}^{\chi},\ZZ_p[G_{\tau}]) \lra \bigwedge^{r-1}\textup{Hom}_{\ZZ_p[\Delta_{\tau^{\prime}}]}(V_{L(\tau^{\prime})}^{\chi},\ZZ_p[G_{\tau^{\prime}}])$$ Therefore $\Phi_0=\{\phi_{\tau}\}_{_{\tau}}$ may be regarded as an element of the group $\hhhk$.
\begin{prop}
\label{homs-hli-1}
Let $\{\phi_{\tau}\}_{_{\tau}}=\Phi_0$ be as above. Then $\phi_{\tau}$ maps $\wedge^r V_{L(\tau)}^{\chi}$ onto $\al_{\tau}^{\chi}(=\textup{ker}(\Psi_{\tau}))$,  for all $\tau$. In particular $\Phi_0$ satisfies $\hli$.
\end{prop}

\begin{proof}
The proof is identical to the proof of Lemma 4.1 of~\cite{kbbonline}, which also follows the proof of Lemma B.1 of~\cite{mr02} line by line.
\end{proof}

\begin{cor}
\label{homs-hli}
Let $\phi$ be an element of $\bigwedge^{r-1}\textup{Hom}_{\ZZ_p[\Delta]}(V_L^{\chi},\ZZ_p)$ such that $$\phi(\wedge^r V_L^{\chi})\subset \al^{\chi}$$ Then there is a $$\Phi^{(0)}=\{\phi_{\tau}^{(0)}\}_{_{\tau}} \in \hhhk$$ which satisfies $\hli$ and is such that $$\phi^{(0)}_1 \mid_{\wedge^r V_L^{\chi}}=\phi \mid_{\wedge^r V_L^{\chi}}$$
\end{cor}

\begin{proof}
If $\phi$ is identically $zero$ the assertion is obvious (choosing $\Phi^{(0)}$ identically $zero$ as well). Otherwise
let $$\{\tilde{\phi}_{\tau}^{(0)}\}_{\tau}=\tilde{\Phi}^{(0)}=[\al^{\chi}:\phi(V_L^{\chi})] \Phi_0$$ where $\Phi_0$ is as above. Then by Proposition~\ref{homs-hli-1}, and noting that the free-of-rank-\emph{one} $\ZZ_p$-module $\al^{\chi}$ has a unique $\ZZ_p$-submodule of index $[\al^{\chi}:~\phi(V_L^{\chi})]$ it follows that $$\tilde{\phi}_1^{(0)}(\wedge^r V_L^{\chi})=\phi(\wedge^r V_L^{\chi})$$ Since $\wedge^r V_L^{\chi}$ is free of rank $one$, we may choose an $\alpha \in \ZZ_p^{\times}$ so that $$\alpha \tilde{\phi}_1^{(0)} \mid_{\wedge^r V_L^{\chi}} \equiv \phi \mid_{\wedge^r V_L^{\chi}}$$ Let now $\Phi^{(0)}:= \alpha \tilde{\Phi}^{(0)}$. Then $\Phi^{(0)} \in \hhhk$ clearly has the desired properties.
\end{proof}


\subsection{Kolyvagin Systems for $(T, \FF_{\al})$:}
\label{kolsys}

In this section, we construct  a Kolyvagin system\footnote{See Definitions 3.1.3 and 3.1.6 of~\cite{mr02} for a more precise definition of a Kolyvagin system.} $\{\kappa_{\tau}^{\Phi}\} \in \textup{\textbf{KS}}(T,\FF_{\al})$ using the Euler system $\{\varepsilon_{L(\tau),\Phi}^{\chi}\}_{\tau}$ of  $\Phi$-Stark elements, for each $\Phi \in \hhhk$ that satisfies $\hli$. We will use these classes in the next section to prove our main results.

Let $\PP$ denote the set primes of $k$ whose elements do not divide $pf_{\chi}$. For each positive integer $m$, let $$\PP_m=\{\qq \in \PP: \qq \hbox{ splits completely in } L(\mu_{p^m})/k\}$$ be a subset of $\PP$. Note that $\PP_m$ is exactly the set of primes being determined by Definition 3.1.6 of~\cite{mr02} or \S 4, Definition 1.1 of~\cite{r00} when $T=\ZZ_p(1)\otimes\chi^{-1}$. Let $\NN$ (\emph{resp.} $\NN_m$) denote the square free products of primes $\qq$ in $\PP$ (\emph{resp.} $\PP_m$).

Let $\FF_{\textup{can}}$ denote the \emph{canonical Selmer structure} given by Definition~ 3.2.1 of \cite{mr02}. Theorem 3.2.4 of~\cite{mr02} gives a map $$\textup{\textbf{ES}}(T)\lra \overline{\textup{\textbf{KS}}}(T,\FF_{\textup{can}},\PP)$$
 where $\textup{\textbf{ES}}(T)$ denotes the collection of Euler systems for $T$, and 
 $$\overline{\textup{\textbf{KS}}}(T,\FF_{\textup{can}},\PP):=\lim_{\stackrel{\longleftarrow}{m}}(\lim_{\stackrel{\lra}{j}} \textup{\textbf{KS}}(T/p^mT,\FF_{\textup{can}},\PP \cap \PP_j))$$ such that if $\{c_{k(\tau)}\}_{\tau}=\textbf{c} \in \textup{\textbf{ES}}(T)$ maps to $\kappa=\left\{\{\kappa_{\tau,m}\}_{\tau \in \NN_{m}}\right\}_m$ under this map (where $\{\kappa_{\tau,m}\}_{\tau \in \NN_{m}} \in \varinjlim_j \textup{\textbf{KS}}(T/p^mT,\FF_{\textup{can}},\PP \cap \PP_j)$ is a Kolyvagin system for $T/p^mT$), then \begin{equation}\label{eqn:ESKS}H^1(k,T)=\varprojlim_{m}H^1(k,T/p^mT) \ni \varprojlim_m \kappa_{1,m}:=\kappa_1=\textbf{c}_{k} \in H^1(k,T).\end{equation}

Let $\boldsymbol{\kappa}^\Phi=\left\{\{\kappa^{\Phi}_{\tau,m}\}_{_{\tau \in \NN_m}}\right\}_m$ be the image of the Euler system $\{\varepsilon_{L(\tau),\Phi}^{\chi}\}_{\tau}$, and let $\kappa_1^{\phi}=\varprojlim \kappa_{1,m}^{\Phi}$ as in (\ref{eqn:ESKS}). In this particular case, (\ref{eqn:ESKS}) reads $$\kappa_1^{\Phi}=\varepsilon_{L,\Phi}^{\chi}.$$

\begin{thm}
\label{mainkolsys}
Assume that $\Phi \in \hhhk$ satisfies $\hli$. Then $$\boldsymbol{\kappa}^\Phi \in \overline{\textup{\textbf{KS}}}(T,\FF_{\al},\PP)$$
\end{thm}

When there is no danger of confusion, we will simply denote the element $\kappa^{\Phi}_{\tau,m} \in L^{\times}/(L^{\times})^{p^m}$ by $\kappa^{\Phi}_{\tau}$. Note that the statement of Theorem~\ref{mainkolsys} says that for each $\tau \in \NN_m$, $\kappa^{\Phi}_{\tau} \in H^1_{\FF_{\al}(\tau)}(k,T/p^mT)$, where $\FF_{\al}(\tau)$ is defined as in Example 2.1.8 of~\cite{mr02}. However, Theorem 3.2.4 of~\cite{mr02} already says that $\kappa^{\Phi}_{\tau} \in  H^1_{\FF_{\textup{can}}(\tau)}(k,T/p^mT)$; therefore to prove Theorem~\ref{mainkolsys} it suffices to prove the following:

\begin{prop}
\label{loc-p}
Let $$\textup{loc}_p: H^1(k,T/p^mT) \lra H^1(k_p,T/p^mT):= \bigoplus_{\wp \mid p}  H^1(k_{\wp},T/p^mT)$$ be the localization map into the semi-local cohomology at $p$. Then $$\textup{loc}_p (\kappa^{\Phi}_{\tau}) \in \al^{\chi}/p^m\al^{\chi} \subset H^1(k_p,T/p^mT)$$
\end{prop}
We prove Proposition~\ref{loc-p} below. We remark that $\al^{\chi}/p^m\al^{\chi}$ is the propagation $H^1_{\FF_{\al}(\tau)}(k_p,T/p^mT)$ of the local condition $H^1_{\FF_{\al}(\tau)}(k_p,T)=\al^{\chi}$ at $p$. Let $$\{\tilde{\kappa}^{\Phi}_{\tau,m} \in H^1(k,T/p^mT) = (L^{\times}/(L^{\times})^{p^m})^{\chi}\}_ {\tau \in \NN_m}$$ be the collection that Definition 4.4.10 of~\cite{r00} associates to the Euler system $\{\varepsilon_{L(\tau),\Phi}^{\chi}\}_{\tau}$. Here we write $\tilde{\kappa}^{\Phi}_{\tau,m}$ for the class denoted by $\kappa_{[k,\tau,m]}$ in~\cite{r00}. We will also often drop $m$ and denote $\tilde{\kappa}^{\Phi}_{\tau,m}$ by $\tilde{\kappa}^{\Phi}_{\tau}$ if there is no danger for confusion. Note that Equation (33) in Appendix A of~\cite{mr02} relates these classes to $\kappa^{\Phi}_{\tau}$.
\begin{lemma}
\label{reduction-1}
Assume $\textup{loc}_p(\tilde{\kappa}^{\Phi}_{\tau}) \in \al^{\chi}/p^m \al^{\chi}$ then $\textup{loc}_p(\kappa^{\Phi}_{\tau}) \in \al^{\chi}/p^m \al^{\chi}$ as well.
\end{lemma}
\begin{proof}
Obvious using Equation (33) in Appendix A of~\cite{mr02}.
\end{proof}

Let $D_{\tau}$ denote the derivative operators as in Definition 4.4.1 of~\cite{r00}. Definition 4.4.10 (and Remark 4.4.3) defines $\tilde{\kappa}^{\Phi}_{\tau}$ as the inverse image of $D_{\tau}\varepsilon_{L(\tau),\Phi}^{\chi}$ (mod~$p^m$) under the restriction map\footnote{ Note that $(\mu_{p^{\infty}}\otimes \chi^{-1})^{G_{k(\tau)}}$ is trivial (where $G_{k(\tau)}$ stands for the absolute Galois group of the totally real field $k(\tau)$),   since, for example, complex conjugation cannot act by $\chi$ on $\mu_{p^{\infty}}$ since $\chi$ is even. This argument proves that this restriction map is an isomorphism, by Remark 4.4.3 of~\cite{r00}.}
\begin{align*} (L^{\times}/(L^{\times})^{p^m})^{\chi}=H^1(k,T/p^mT) &\lra H^1(k(\tau),T/p^mT)^{G_{\tau}}\\
&=[(L(\tau)^{\times}/(L(\tau)^{\times})^{p^m})^{\chi}]^{G_{\tau}}
\end{align*} 
Therefore $\textup{loc}_p(\tilde{\kappa}^{\Phi}_{\tau})$ maps to $\textup{loc}_p(D_{\tau}\varepsilon_{L(\tau),\Phi}^{\chi})$ (mod~$p^m$) under the map (which is also an isomorphism by Remark 4.4.3, Proposition B.5.1 and Proposition B.4.2 of~\cite{r00}) 
$$H^1(k_p,T/p^mT) \lra H^1(k(\tau)_p,T/p^mT)^{G_{\tau}}$$ 
Under this isomorphism $V_L^{\chi}/p^mV_L^{\chi}$ is mapped isomorphically onto \\$(V_{L(\tau)}^{\chi}/p^mV_{L(\tau)}^{\chi})^{G_{\tau}}$ and hence, by definition of $\al_{\tau}^{\chi}$ and the fact that it is a free $\ZZ_p[G_{\tau}]$-module, $\al^{\chi}/p^m\al^{\chi}$ is mapped isomorphically onto $[\al_{\tau}^{\chi}/p^m\al_{\tau}^{\chi}]^{G_{\tau}}$. The diagram below summarizes the discussion in this paragraph:
$$\xymatrix{
H^1(k_p,T/p^mT) \ar[r]^(.45){\sim}& H^1(k(\tau)_p,T/p^mT)^{G_{\tau}}\\
V_L^{\chi}/p^mV_L^{\chi} \ar[r]^(.42){\sim} \ar@{^{(}->}[u] &(V_{L(\tau)}^{\chi}/p^mV_{L(\tau)}^{\chi})^{G_{\tau}} \ar@{^{(}->}[u]\\
\al^{\chi}/p^m\al^{\chi}\ar[r]^(.4){\sim} \ar@{^{(}->}[u] &[\al_{\tau}^{\chi}/p^m\al_{\tau}^{\chi}]^{G_{\tau}}\ar@{^{(}->}[u]
}$$

\begin{prop}
\label{lined}
Suppose $\Phi$ satisfies $\hli$. Then $$\textup{loc}_p(\tilde{\kappa}^{\Phi}_{\tau}) \in \al^{\chi}/p^m \al^{\chi}.$$
\end{prop}

\begin{proof}
Since $\textup{loc}_p$ is Galois equivariant $\textup{loc}_p(D_{\tau}\varepsilon_{L(\tau),\Phi}^{\chi})=D_{\tau}\textup{loc}_p(\varepsilon_{L(\tau),\Phi}^{\chi})$. Further $\textup{loc}_p(\varepsilon_{L(\tau),\Phi}^{\chi}) \in \al_{\tau}^{\chi}$ since $\Phi$ satisfies $\hli$. On the other hand, by Lemma 4.4.2 of~\cite{r00} $D_{\tau}\varepsilon_{L(\tau),\Phi}^{\chi}$ (mod~$p^m$) is fixed by $G_{\tau}$, which in return implies $$\textup{loc}_p(\varepsilon_{L(\tau),\Phi}^{\chi}) \,(\textup{mod\,} p^m) \in [\al_{\tau}^{\chi}/p^m\al_{\tau}^{\chi}]^{G_{\tau}}$$ This proves that $\textup{loc}_p(\tilde{\kappa}^{\Phi}_{\tau})$ maps into $\al^{\chi}/p^m \al^{\chi}$ by above discussion.
\end{proof}
\begin{proof}[Proof of Proposition~\ref{loc-p}]
Immediately follows from Lemma~\ref{reduction-1} and Proposition~\ref{lined}.
\end{proof}
By the discussion following the statement of Theorem~\ref{mainkolsys}, this also completes the proof of Theorem~\ref{mainkolsys}.
\newpage
\section{Applications to Ideal Class Groups}
\label{sec:applications}

In this section, we apply the Kolyvagin system machinery developed in~\cite{mr02} to the Kolyvagin system we constructed in~\S\ref{kolsys}. Let $$\iota: (O_L^\times)^{\chi} \lra V_L^{\chi}$$ be the localization map at $p$. The induced map $\wedge^r (O_L^\times)^{\chi} \ra \wedge^r V_L^{\chi}$ will be denoted by $\iota^{(r)}$. Throughout this section we assume the following hypotheses:
\begin{itemize}
\item [\textbf{H-S}] The finite set $S$ contains no finite primes which split completely in $L/k$.
    
\item [\textbf{H-F}] $[\wedge^r V_L^{\chi} : \ZZ_p\cdot\iota^{(r)}(\varepsilon_{L}^{\chi})]<\infty$.
\end{itemize}

Since $\varepsilon_{L}^{\chi} \in \wedge^{r}(O_{L}^{\times})^{\chi}$ by Remark~\ref{integrality}, hypothesis \textbf{H-F} makes sense.
\begin{rem}
\label{rem:finite-index}
The proof of Proposition 6.6 of~\cite{ru96} shows that \textbf{H-F} holds  if one assumes \textbf{H-S} and that the map $\iota$ above is injective. Note that Leopoldt's conjecture (\textbf{LC}) for $L$ guarantees that $\iota$ is injective, thus implies \textbf{H-F}. 
\end{rem}

Recall that for an appropriately chosen (namely those satisfying $\hli$) $$\Phi \in \varprojlim_{L(\tau) \in \kk}\bigwedge^{r-1} \hbox{Hom}_{\ZZ_p[\Delta_{\tau}]}(V_{L(\tau)}^{\chi},\ZZ_p[\Delta_\tau]^{\chi} )$$ we obtained a Kolyvagin system $\{\kappa_{\tau}^{\Phi}\}$ for the Selmer structure $(T, \FF_{\al}, \PP)$.

The following theorem (which is really Corollary 5.2.13 of~\cite{mr02}) is the main application of the Kolyvagin systems we constructed. Notice that Corollary 5.2.13 of~\cite{mr02} is valid by Proposition~\ref{modifiedcorerank}.

\begin{thm}
\label{bound}
$$\textup{length}_{\ZZ_p} (H^1_{\FF^*_{\al}}(k,T^*)) \leq \max\{j: \kappa_1^{\Phi} \in p^j H^1_{\FF_{\al}}(k,T)\}$$ for all $\Phi$ satisfying $\hli$. We have equality if the Kolyvagin system $\{\kappa_{\tau}^{\Phi}\}$ is primitive in the sense of Definition 3.4.5 of~\cite{mr02}.
\end{thm}

Using the computations in Section~\ref{sec:kolsys} we prove



\begin{prop}
\label{index-kolsys} Let $\tau,\Delta_{\tau}, V_K$ and $\iota$ be defined as above. Then
\begin{align*}
\{\iota(\varepsilon_{L,\Phi}^{\chi}): \Phi \in \varprojlim_\tau \bigwedge^{r-1} \textup{Hom}_{\ZZ_p[\Delta_{\tau}]}(V_{L(\tau)}^{\chi}, \ZZ_p[&\Delta_{\tau}]^{\chi})\}\\&=[\wedge^r V_L^{\chi} : \ZZ_p\cdot\iota^{(r)}(\varepsilon_{L}^{\chi})]V_L^{\chi}.
\end{align*}
\end{prop}

\begin{proof}
By Proposition~\ref{liftinghoms} and Lemma~\ref{krasner} $$\{\Phi_1(\mathbf{v}): \mathbf{v} \in \wedge^r V_L^{\chi}, \{\Phi_{\tau}\}_{\tau} \in \varprojlim_{\tau} \bigwedge^{r-1}\textup{Hom}_{\ZZ_p[\Delta_{\tau}]}(V_{L(\tau)}^{\chi}, \ZZ_p[\Delta_{\tau}]^{\chi})\}=V_L^{\chi}.$$ (See also the proof of~\cite{ru96} Proposition 6.6.) On the other hand, $V_L$ is free of rank $r$ over $\ZZ_p[\hbox{Gal}(L/k)]$ by Krasner's lemma, therefore $V_L^{\chi}$ is free of rank $r$ over $\ZZ_p$. This implies that $\wedge^r V_L^{\chi}$ is free of rank \emph{one} over $\ZZ_p$. The proof of the proposition now follows.
\end{proof}
\begin{lemma}
\label{lemma-L-component}
\begin{align*} &\{\iota(\varepsilon_{L,\Phi}^{\chi}): \Phi \in \varprojlim_{\tau} \bigwedge^{r-1} \textup{Hom}_{\ZZ_p[\Delta_{\tau}]}(V_{L(\tau)}^{\chi}, \ZZ_p[\Delta_{\tau}]^{\chi}), \Phi \hbox{ satisfies } \hli\}=\\ &=\{\iota(\varepsilon_{L,\Phi}^{\chi}): \Phi \in \varprojlim_{\tau} \bigwedge^{r-1} \textup{Hom}_{\ZZ_p[\Delta_{\tau}]}(V_{L(\tau)}^{\chi}, \ZZ_p[\Delta_{\tau}]^{\chi})\} \cap \al^{\chi}
\end{align*}
\end{lemma}
\begin{proof}
The left hand side is clearly contained in the right hand side. Let $$c \in \{\iota(\varepsilon_{L,\Phi}^{\chi}): \Phi \in \varprojlim_{\tau} \bigwedge^{r-1} \textup{Hom}_{\ZZ_p[\Delta_{\tau}]}(V_{L(\tau)}^{\chi}, \ZZ_p[\Delta_{\tau}]^{\chi})\} \cap \al^{\chi}$$ be an element of the set on the right hand side, so that $$c=\phi\circ\iota(\varepsilon_{L,\Phi}^{\chi}) \in \al^{\chi}$$ for some $$\phi \in \wedge^{r-1}\textup{Hom}_{\ZZ_p[\Delta]}(V_L^{\chi},\ZZ_p[\Delta]^{\chi}).$$ Since we assumed \textbf{H-F}, $\ZZ_p\cdot\iota^{(r)}(\varepsilon_{L}^{\chi})$ is a subgroup of finite index in $\wedge^r V_L^{\chi}$ and since $V_L^{\chi}/\al^{\chi}$ is torsion free, it follows that $\phi(\wedge^r V_L^{\chi}) \subset \al^{\chi}$. Then by Corollary~\ref{homs-hli} there is a $$\begin{array}{rcl}\{\phi_{\tau}^{(0)}\}_{\tau}=\Phi^{(0)} &\in& \hhhk\\ &=&  \varprojlim_\tau \bigwedge^{r-1} \textup{Hom}_{\ZZ_p[\Delta_{\tau}]}(V_{L(\tau)}^{\chi}, \ZZ_p[\Delta_{\tau}]^{\chi})\end{array}$$ which satisfies $\hli$ and such that $$\phi_1^{(0)} \mid_{\wedge^r V_L^{\chi}} \equiv\phi \mid_{\wedge^r V_L^{\chi}}$$ This shows $c=\varepsilon_{L,\Phi^{(0)}}^{\chi}$ and thus $c$ belongs to the left hand side as well.

\end{proof}

\begin{cor}
\label{i-component}
\begin{align*} \{\iota(\varepsilon_{L,\Phi}^{\chi}): \Phi \in \varprojlim_{\tau} \bigwedge^{r-1} \textup{Hom}_{\ZZ_p[\Delta_{\tau}]}(V_{L(\tau)}^{\chi}, \ZZ_p&[\Delta_{\tau}]^{\chi}), \Phi \hbox{ satisfies } \hli\}\\ &=[\wedge^r V_L^{\chi} : \ZZ_p\cdot\iota^{(r)}(\varepsilon_{L}^{\chi})]\al^{\chi}\\
\end{align*}
\end{cor}

\begin{proof}By Lemma~\ref{lemma-L-component} the left hand side above equals $$\{\iota(\varepsilon_{L,\Phi}^{\chi}): \Phi \in \varprojlim_{\tau} \bigwedge^{r-1} \textup{Hom}_{\ZZ_p[\Delta_{\tau}]}(V_{L(\tau)}^{\chi}, \ZZ_p[\Delta_{\tau}]^{\chi})\} \cap \al^{\chi}
$$ which, by Proposition~\ref{index-kolsys}, is equal to $\left\{[\wedge^r V_L^{\chi} : \ZZ_p\cdot\iota^{(r)}(\varepsilon_{L}^{\chi})]V_L^{\chi}\right\} \cap \al^{\chi}$, and this equals $[\wedge^r V_L^{\chi} : \ZZ_p\cdot\iota^{(r)}(\varepsilon_{L}^{\chi})]\al^{\chi}$ since $V_L^{\chi}/\al^{\chi}$ is $\ZZ_p$-torsion-free.

\end{proof}

Note that Corollary~\ref{i-component} implies in particular that $\varepsilon_{L,\Phi}^{\chi}\neq0$ for some choice of $\Phi$ which satisfies $\hli$. This, together with Theorem~\ref{bound} and the right-most equality in (\ref{eqn:ESKS}) proves

\begin{cor}
\label{finiteness}
$H^1_{\FF_{\al}^{*}}(k,T^*)$ is finite.
\end{cor}
\begin{rem}
There is of course a direct way to prove Corollary~\ref{finiteness} using Remark~\ref{rem:LCshows} and the finiteness of the ideal class group, if one assumes \textbf{LC} for $L$.
\end{rem}

Until the end of this section we assume \textbf{LC} for the number field $L$. This in particular means that the canonical map $$\iota:(O_L^{\times})^{\chi} \lra V_L^{\chi}$$ is injective. We will then identify $(O_L^{\times})^{\chi}$ with its image inside $V_L^{\chi}$ under the map $\iota$, and drop $\iota$ from notation.

Recall the exact sequence
$$0 \lra (O_L^{\times})^{\chi}\cap \al^{\chi} \lra (O_L^{\times})^{\chi} \lra V_L^{\chi}/{\al^{\chi}}. $$ We note that the intersection $(O_L^{\times})^{\chi}\cap \al^{\chi}$ is taken inside $V_L^\chi$ after identifying $(O_L^{\times})^{\chi}$ with its image inside $V_L^{\chi}$, and is free of rank \emph{one} (by rank considerations in the above sequence).

\begin{thm}
\label{trueindexformodified}
\begin{align*} \{\varepsilon_{L,\Phi}^{\chi}: \Phi \in \varprojlim_{\tau} \bigwedge^{r-1} \textup{Hom}_{\ZZ_p[\Delta_{\tau}]}(V_{L(\tau)}^{\chi}, &\ZZ_p[\Delta_{\tau}]^{\chi}), \,\Phi \hbox{ satisfies } \hli\}=\\ &=\frac{[\wedge^r V_L^{\chi} : \ZZ_p\varepsilon_{L}^{\chi}]}{[\al^{\chi}:(O_L^{\times})^{\chi}\cap \al^{\chi}]}(O_L^{\times})^{\chi}\cap \al^{\chi}.\\
\end{align*}

\textup{(}Note that we still use the additive notation for the multiplicative group $(O_L^{\times})^{\chi}$\textup{)}.
\end{thm}

\begin{proof}
This follows from Corollary~\ref{i-component} together with  our assumption \textbf{LC} and that $\al^{\chi}$ is a free $\ZZ_p$-module of rank \emph{one}.

\end{proof}

\begin{cor}
\label{application-modified}
$$|H^1_{\FF^*_{\al}}(k,T^*)| \leq \frac{[\wedge^r V_L^{\chi} : \ZZ_p\varepsilon_{L}^{\chi}]}{[\al^{\chi}:(O_L^{\times})^{\chi}\cap \al^{\chi}]}.$$

\end{cor}

\begin{proof}
By Theorem~\ref{bound} it suffices to show that $$\min_{\{\Phi \textup{ satisfying } \tiny\hli \}} \max\{p^j: \varepsilon_{L,\Phi}^{\chi} \in p^j H^1_{\FF_{\al}}(k,T)\}=\frac{[\wedge^r V_L^{\chi} : \ZZ_p\varepsilon_{L}^{\chi}]}{[\al^{\chi}:(O_L^{\times})^{\chi}\cap \al^{\chi}]}.$$ This is exactly what Theorem~\ref{trueindexformodified} says, recalling that $$H^1_{\FF_{\al}}(k,T)=(O_L^{\times})^{\chi}\cap \al^{\chi}$$ by Proposition~\ref{modifiedexplicit}(i).
\end{proof}

We now use Corollary~\ref{application-modified}, Proposition~\ref{comparison} and Remark~\ref{rem:LCshows} to obtain bounds on the classical Selmer group.

\begin{thm}
\label{main}
$|A_L^{\chi}| \leq [\wedge^r (O_L^{\times})^{\chi}: \ZZ_p\varepsilon_{L}^{\chi}]$
\end{thm}

\begin{proof}
It follows from  Proposition~\ref{explicit}, Corollary~\ref{application-modified} and Remark~\ref{rem:LCshows} that $$\mid A_L^{\chi}\mid \leq \frac{[\wedge^r V_L^{\chi} : \ZZ_p\varepsilon_{L}^{\chi}]}{[V_L^{\chi}:(O_L^{\times})^{\chi}]}$$
On the other hand, since $V_L^{\chi}$ and the image of $ (O_L^{\times})^{\chi}$ inside $V_L^{\chi}$ are both free of rank $r$, it follows that  $$[\wedge^r V_L^{\chi}:\wedge^r (O_L^{\times})^{\chi}]=[V_L^{\chi}:(O_L^{\times})^{\chi}]$$
 But this means $$\frac{[\wedge^r V_L^{\chi} : \ZZ_p\varepsilon_{L}^{\chi}]}{[V_L^{\chi}:(O_L^{\times})^{\chi}]}=\frac{[\wedge^r V_L^{\chi} : \ZZ_p\varepsilon_{L}^{\chi}]}{[\wedge^r V_L^{\chi} : \wedge^r (O_L^{\times})^{\chi}]}=[\wedge^r (O_L^{\times})^{\chi}: \ZZ_p\varepsilon_{L}^{\chi}]$$ and the proof of the theorem follows.
\end{proof}
Recall our assumption on the set $S$ that it does not contain any non-archimedean prime of $k$ which splits completely in $L/k$. Fix also an embedding  of $\QQ_p$ into the complex numbers $\mathbb{C}$. For $x,y \in \mathbb{C}^{\times}$ we will write $x \sim y$ if $x/y \in \ZZ_p^{\times}$.

By their conjectural description, the Stark elements of Rubin relate to the $L$-values:
\begin{equation}
\label{eqn:starkL}
R_{\chi}(\varepsilon_{L}^{\chi})\sim \lim_{s \ra 0}s^{-r}L_{S_{\infty},T}(s,\chi^{-1})
\end{equation}
where $R_{\chi}$ is a regulator map (as described in \S 2.1 of~\cite{ru96}) on $(\wedge^r O_{L,S_{\infty},T}^{\times})^{\chi}=(\wedge^r O_{L,T}^{\times})^{\chi}$, $S_{\infty}$ is the set of infinite places of $L$ and $L_{S_{\infty},T}(s,\chi^{-1})$ is the modified Artin $L$-function attached to $\chi^{-1}$, and $T$ is a set of primes discussed at the beginning of Section~\ref{subsec:stark}.  Note that we use $L_{S_{\infty},T}(s,\chi^{-1})$ instead of the modified Artin $L$-function $L_{S,T}(s,\chi^{-1})$, in the expense of replacing the equality in~\ref{eqn:starkL} by "$\sim$". We further remark that $$\lim_{s \ra 0}s^{-r}L_{S,T}(s,\chi^{-1}) \sim \lim_{s \ra 0}s^{-r}L_{S_{\infty},T}(s,\chi^{-1})$$ since we insist that $S$ does not contain any non-archimedean prime of $k$ which splits completely in $L/k$.

On the other hand we have the analytic class number formula (see \cite{gross})
\begin{equation}
\label{analytic}
\lim_{s \ra 0} s^{1-\#(S_{\infty})} \zeta_{L,S_{\infty},T}(s)=-\frac{\#(A_{S_{\infty},T})R_{S_{\infty},T}}{\#(\mu_T)}
\end{equation}
 See  \S1.1 of~\cite{ru96} for definitions of $A_{S_{\infty},T}, R_{S_{\infty},T}$ and $\mu_T$. As we have noted at the beginning of Section~\ref{subsec:stark}, one may choose the set of primes $T$ in a way that $\#\mu_T =1$ and also that $$\lim_{s \ra 0} s^{1-\#(S_{\infty})} \zeta_{L,S_{\infty},T}(s) \sim \lim_{s \ra 0} s^{1-\#(S_{\infty})} \zeta_{L,S_{\infty}}(s),$$ $$\#(A_{S_{\infty},T}) \sim \#(A_{S_{\infty}}), R_{S_{\infty},T} \sim R_{S_{\infty}} $$
Using~(\ref{eqn:starkL}) for all non-trivial characters of $\Delta$ and (\ref{analytic}) together with above  remarks and Theorem~\ref{main} yields (see \S 5 of~\cite{ru92} for more details):

\begin{thm}
\label{thm:mainequality}
$|A_L^{\chi}| = [\wedge^r (O_L^{\times})^{\chi}: \ZZ_p\varepsilon_{L}^{\chi}]$
\end{thm}
This in particular shows that for $\Phi_0$ as defined in~\S\ref{homs} the Kolyvagin system $\{\kappa_{\tau}^{\Phi_0}\}$ is $primitive$ if the set $S$ contains no non-archimedean prime of $k$ which splits completely in $L/k$.

Let $\mathcal{R}_L^\chi$ be the value of the regulator map $R_\chi$ evaluated at a generator of $\wedge^r(O_L^\times)^{\chi}$. Theorem~\ref{thm:mainequality} together with (\ref{eqn:starkL}) implies the following class number formula:
\begin{cor}
$$\lim_{s\ra0} s^{-r}L_{S_{\infty},T}(s,\chi^{-1}) \sim \mathcal{R}_L^\chi\cdot |A_L^\chi|.$$
\end{cor}
\bibliographystyle{plain}
\bibliography{references}

\begin{thebibliography}{10}

\bibitem{kbbonline}
K\^az{\i}m B{\"u}y\"ukboduk.
\newblock {$\Lambda$}-adic {K}olyvagin systems, 2007.
\newblock 56pp., preprint. \textup{http://arxiv.org/abs/0706.0377v1}.

\bibitem{gross}
Benedict~H. Gross.
\newblock On the values of abelian {$L$}-functions at {$s=0$}.
\newblock {\em J. Fac. Sci. Univ. Tokyo Sect. IA Math.}, 35(1):177--197, 1988.

\bibitem{krasner}
Marc Krasner.
\newblock Sur la repr\'esentation exponentielle dans les corps relativement
  galoisiens de nombers $\hbox{{\frakfamily p}}$-adiques.
\newblock {\em Acta Arith.}, 3:133--173, 1939.

\bibitem{mr02}
Barry Mazur and Karl Rubin.
\newblock Kolyvagin systems.
\newblock {\em Mem. Amer. Math. Soc.}, 168(799):viii+96, 2004.

\bibitem{pr-es}
Bernadette Perrin-Riou.
\newblock Syst\`emes d'{E}uler {$p$}-adiques et th\'eorie d'{I}wasawa.
\newblock {\em Ann. Inst. Fourier (Grenoble)}, 48(5):1231--1307, 1998.

\bibitem{popescu}
Cristian~D. Popescu.
\newblock Gras-type conjectures for function fields.
\newblock {\em Compositio Math.}, 118(3):263--290, 1999.

\bibitem{ru92}
Karl Rubin.
\newblock Stark units and {K}olyvagin's ``{E}uler systems''.
\newblock {\em J. Reine Angew. Math.}, 425:141--154, 1992.

\bibitem{ru96}
Karl Rubin.
\newblock A {S}tark conjecture ``over {$\mathbf Z$}'' for abelian
  {$L$}-functions with multiple zeros.
\newblock {\em Ann. Inst. Fourier (Grenoble)}, 46(1):33--62, 1996.

\bibitem{r00}
Karl Rubin.
\newblock {\em Euler systems}, volume 147 of {\em Annals of Mathematics
  Studies}.
\newblock Princeton University Press, Princeton, NJ, 2000.
\newblock Hermann Weyl Lectures. The Institute for Advanced Study.

\bibitem{solomon-2}
David Solomon.
\newblock On {$p$}-adic abelian {S}tark conjectures at {$s=1$}.
\newblock {\em Ann. Inst. Fourier (Grenoble)}, 52(2):379--417, 2002.

\bibitem{solomon-1}
David Solomon.
\newblock Twisted zeta-functions and abelian {S}tark conjectures.
\newblock {\em J. Number Theory}, 94(1):10--48, 2002.

\bibitem{tate}
John Tate.
\newblock {\em Les conjectures de {S}tark sur les fonctions {$L$} d'{A}rtin en
  {$s=0$}}, volume~47 of {\em Progress in Mathematics}.
\newblock Birkh\"auser Boston Inc., Boston, MA, 1984.
\newblock Lecture notes edited by Dominique Bernardi and Norbert Schappacher.

\bibitem{wiles}
Andrew Wiles.
\newblock Modular elliptic curves and {F}ermat's last theorem.
\newblock {\em Ann. of Math. (2)}, 141(3):443--551, 1995.

\end{thebibliography}
\end{document}